\documentclass{amsart}
\usepackage{mathbbol}
\usepackage{amsmath, amsthm, amssymb}

\newtheorem{thm}{Theorem}[section]
\newcommand{\bt}{\begin{thm}}
\newcommand{\et}{\end{thm}}

\newtheorem{cor}[thm]{Corollary}
\newcommand{\bc}{\begin{cor}}
\newcommand{\ec}{\end{cor}}

\newtheorem{lem}[thm]{Lemma}
\newcommand{\bl}{\begin{lem}}
\newcommand{\el}{\end{lem}}

\newtheorem{prop}[thm]{Proposition}
\newcommand{\bp}{\begin{prop}}
\newcommand{\ep}{\end{prop}}

\newtheorem{defn}[thm]{Definition}
\newcommand{\bd}{\begin{defn}}      
\newcommand{\ed}{\end{defn}}

\newtheorem{rmrk}[thm]{Remark}
\newcommand{\br}{\begin{rmrk}}
\newcommand{\er}{\end{rmrk}}

\newtheorem{example}[thm]{Example}

\newcommand{\thmref}[1]{Theorem~\ref{#1}}
\newcommand{\secref}[1]{Section~\ref{#1}}
\newcommand{\lemref}[1]{Lemma~\ref{#1}}
\newcommand{\defref}[1]{Definition~\ref{#1}}

\newcommand{\propref}[1]{Proposition~\ref{#1}}

\newcommand{\exref}[1]{Example~\ref{#1}}

\newcommand{\N}{\mathbb{N}}

\newcommand{\R}{\mathbb{R}}
\newcommand{\Z}{\mathbb{Z}}
\newcommand{\dist}{\operatorname{dist}}
\newcommand{\diam}{\operatorname{diam}}

\newcommand{\hm}{{\mathcal H}}

\newcommand{\lip}{\operatorname{Lip}}
\newcommand{\mass}[2][]{{\mathbf M_{#1}}(#2)}

\newcommand{\form}{{\mathcal D}}        
\newcommand{\curr}{{\mathbf M}}         
\newcommand{\norcurr}{{\mathbf N}}      
\newcommand{\rectcurr}{{\mathcal R}}    
\newcommand{\intrectcurr}{{\mathcal I}} 
\newcommand{\intcurr}{{\mathbf I}}      

\newcommand{\nmass}[1]{{\mathbf N}(#1)}

\newcommand{\iflnorm}[1]{{\mathcal F}(#1)}
\newcommand{\ifillnorm}[1]{{\operatorname{Fillvol}}(#1)}
\newcommand{\fillrad}[1]{{\operatorname{Fillrad}}(#1)}

\newcommand{\rstr}{\:\mbox{\rule{0.1ex}{1.2ex}\rule{1.1ex}{0.1ex}}\:}
\newcommand{\bdry}{\partial\hspace{-0.05cm}}
\newcommand{\slice}[3]{\langle#1,#2,#3\rangle}

\newcommand{\spt}{\operatorname{spt}}
\newcommand{\ohne}{\backslash}

\begin{document}

\title[Flat convergence for metric integral currents]{Flat convergence for integral currents in metric spaces}

\author{Stefan Wenger}

\address
  {Department of Mathematics,
   University of Basel,
   Rheinsprung 21,
   CH -- 4051 Basel,
   Switzerland}
\email{stefan.wenger@unibas.ch}

\date{February 15, 2005}




\maketitle

\section{Introduction}
It is well known that in compact local Lipschitz neighborhood retracts in $\R^n$ flat convergence for Euclidean integer rectifiable currents amounts just to 
weak convergence. 
The purpose of the present paper is to extend this result to integral currents in complete metric spaces admitting a local cone type inequality. This includes
for example all Banach spaces and complete ${\rm CAT}(\kappa)$-spaces, $\kappa\in\R$. The main result can be used e.g.~to prove the existence of minimal elements
in a fixed Lipschitz homology class in compact metric spaces admitting local cone type inequalities or to conclude that integral currents which are weak limits of 
sequences of absolutely area minimizing integral currents are again absolutely area minimizing.
%
\subsection{Statement of the main results}
The theory of normal and integral currents in Euclidean space was developed by Federer and Fleming in \cite{Fed-Flem}.
Recently, using ideas of De Giorgi \cite{deGiorgi}, Ambrosio and Kirchheim \cite{Ambr-Kirch-curr} extended this theory to the setting of arbitrary 
complete metric spaces. In
the new theory the space of differential $k$-forms, used in the classical theory to define $k$-dimensional currents, is replaced by the space
\begin{equation*}
 \form^k(X):= \{(f,\pi_1,\dots,\pi_k): \text{ $f,\pi_i:X\to \R$ Lipschitz, $f$ bounded}\}.
\end{equation*}
By definition, a $k$-dimensional metric current is a multi-linear functional on $\form^k(X)$ satisfying a continuity, a locality and a finite mass condition. See 
\secref{section:currents} for all definitions from the theory of metric currents relevant for our purposes. The space $\intcurr_k(X)$ of $k$-dimensional metric 
integral currents in $X$ roughly consists of those $k$-dimensional metric currents in $X$ that correspond to $k$-dimensional rectifiable sets with orientation
and multiplicities.

There are two important notions of convergence for (metric) integral currents. A sequence $(T_m)\subset\intcurr_k(X)$ is said to weakly converge to some 
$T\in\intcurr_k(X)$ if
\begin{equation*}
 T_m(f,\pi_1,\dots,\pi_k)\to T(f,\pi_1,\dots,\pi_k)
\end{equation*}
for all $(f,\pi_1,\dots,\pi_k)\in\form^k(X)$, i.e.~if $T_m$ converges pointwise to $T$. A more geometric notion of 
convergence is that with respect to the flat distance given
for $T_1,T_2\in\intcurr_k(X)$ by $d_{\mathcal F}(T_1,T_2):= \iflnorm{T_1-T_2}$ where
\begin{equation*}
  \iflnorm{T}:=\inf\{\mass{U}+\mass{V}:T=U+\bdry V, U\in\intcurr_k(X), V\in\intcurr_{k+1}(X)\}.
\end{equation*}
An important and useful result in the classical theory states that in compact Lipschitz neighborhood retracts in $\R^n$ weak convergence is equivalent to convergence with 
respect to the flat distance. (We will sometimes call the latter flat convergence.)
The proof of this result relies on the deformation theorem which gives a way of deforming a current into the $k$-dimensional skeleton of the standard cubic 
subdivision of $\R^n$ with bounded increase of mass. 
In general, there is no analogue of the deformation theorem in metric spaces. The purpose of this paper is to extend the equivalence of weak and flat convergence
to quasiconvex metric spaces admitting local cone type inequalities. This includes, in particular, infinite dimensional Banach spaces (for which there cannot exist
the deformation theorem from Euclidean space). Other examples of metric spaces satisfying the two conditions will be given later 
(see eg. \secref{section:local-cone-inequalities}).

In connection with the flat distance the notion of integral filling volume will be important. For a $T\in\intcurr_k(X)$ it is defined by
\begin{equation*}
 \ifillnorm{T}:= \inf\{\mass{V}: \text{$V\in\intcurr_{k+1}(X)$ and $\bdry V =  T$}\}
\end{equation*}
where we use the convention that the infimum of the empty set be infinite.
Clearly, we have 
\begin{equation*}
 \iflnorm{T}\leq\min\{\ifillnorm{T},\mass{T}\}
\end{equation*}
for all $T\in\intcurr_k(X)$, $k\geq 0$.
In general the inequality is strict, as easy examples show.

For notational purposes and for convenience in the proof we deal with $0$-dimensional and higher dimensional currents separately.

\bd\label{definition:quasi-convex}
 A metric space $(X,d)$ is said to be $\delta$-quasiconvex if there exists a constant $C<\infty$ such that every two points $x,y\in X$ with $d(x,y)\leq\delta$ 
 can be joined by a Lipschitz curve $\gamma_{xy}:[0,1]\to X$ of length at most $Cd(x,y)$.
\ed
A quasiconvex metric space is a space which is $\delta$-quasiconvex for some $\delta>0$.

\bt\label{theorem:0-dimension-flat-weak}
 Let $(X,d)$ be a complete quasiconvex metric space and $T\in\intcurr_0(X)$. Then a bounded sequence $(T_m)\subset\intcurr_0(X)$ converges weakly to 
 $T$ if and only if 
 \begin{equation*}
  \ifillnorm{T-T_m}\to 0\quad\text{as $m\to\infty$}.
 \end{equation*}
\et
Here, a sequence $(T_m)\subset\intcurr_k(X)$ is said to be bounded if the sequence ($\mass{T_m}$) and, in case $k\geq 1$, also the sequence ($\mass{\bdry T_m})$
is bounded.

As for the higher dimensional analogue of \thmref{theorem:0-dimension-flat-weak} we need the definition of cone type inequality.
\fussy
\bd
 Let $(X,d)$ be a complete metric space, $k\in\N$, $\delta>0$, and $C>0$. Then $X$ is said to admit a $\delta$-cone type inequality for $\intcurr_k(X)$ 
 with constant $C$ if for every $T\in\intcurr_k(X)$ with $\bdry T=0$ and $\diam(\spt T)\leq \delta$ there exists an $S\in\intcurr_{k+1}(X)$ satisfying
 $\bdry S=T$ and
 \begin{equation*}
  \mass{S}\leq C\diam(\spt T)\mass{T}.
 \end{equation*}
\ed
We say that a space $X$ admits a local cone type inequality for $\intcurr_k(X)$ if it admits a $\delta$-cone type inequality  for $\intcurr_k(X)$ 
for some $\delta>0$. If $\delta=\infty$ then we say that $X$ admits a global cone type inequality for $\intcurr_k(X)$. Spaces admitting local cone type inequalities
will be dealt with in \secref{section:local-cone-inequalities}. As mentioned above all Banach spaces and all ${\rm CAT}(\kappa)$-spaces, $\kappa\in\R$, admit local 
cone type inequalities. Furthermore, compact Lipschitz neighborhood retracts in $\R^n$ also admit local cone type inequalities.
\bt\label{theorem:flat-weak-convergence-higher-dimension}
 Let $(X,d)$ be a complete quasiconvex metric space and $k\geq 1$. Suppose that $X$ admits local cone type 
 inequalities for $\intcurr_j(X)$ for $j=1,\dots,k$. If $(T_m)_{m\in\N}\subset\intcurr_k(X)$ is a bounded sequence and 
 $T\in\intcurr_k(X)$ then $T_m$ weakly converges to $T$ if and only if
\begin{equation*}
 \lim_{m\to\infty}\iflnorm{T-T_m}=0.
\end{equation*}
Moreover, if $\bdry T_m=0$ for all $m\in\N$ then $T_m$ weakly converges to $T$ if and only if $\ifillnorm{T-T_m}\to 0$.
\et

We point out that we do not make any further assumptions on the $T_m$. In particular, $T_m$ can have unbounded support.
We also mention that without any compactness assumption, the analogous statement for Euclidean currents 
in $\R^n$ weakly converging to $0$ in the usual sense with compactly supported differential forms is false, as easy examples show.


As mentioned above, Theorems \ref{theorem:0-dimension-flat-weak} and \ref{theorem:flat-weak-convergence-higher-dimension} are of use in connection with 
area minimizing integral currents. We give one prominent application here.
For a complete metric space 
$(X,d)$ and $A\subset X$ closed we denote by ${\bf H}_k(X,A)$ the $k$-th homology class of compactly supported integral currents with boundary in $A$ 
(see \secref{section:applications} for the precise definition).
\bt\label{theorem:minimal-homology}
 Let $k\geq 1$ be an integer and $(X,d)$ a quasiconvex compact metric space and $A\subset X$ closed and quasiconvex (when endowed with the induced metric). 
 Suppose furthermore that both $X$ and $A$ admit local cone type inequalities for $\intcurr_j(X)$ and $\intcurr_j(A)$, 
 for $j=1,\dots,k$. Then, for every homology class $c\in {\bf H}_k(X,A)$ there exists a $T\in{\mathcal Z}_k(X,A)$ with $[T]=c$ and such that
\begin{equation*}
 \mass{T} = \inf\{\mass{T'}:T'\in{\mathcal Z}_k(X,A), [T']= c\}.
\end{equation*}
\et
The proof, as well as further applications, will be given in \secref{section:applications}.
\subsection{Outline of the proof of \thmref{theorem:flat-weak-convergence-higher-dimension}} The proof is by induction on $k$ and the interesting case
is that of $k\geq 2$. We only sketch the proof of the second statement since the first statement readily follows from the second one.
Our strategy will be to find for $\varepsilon>0$ and for each bounded sequence $(T_m)\subset\intcurr_k(X)$ converging weakly to $0$ decompositions
\begin{equation*}
 T_m = U_m^1+\dots + U_m^{M_m}+R_m
\end{equation*}
into the sum of integral cycles such that
\begin{enumerate}
 \item $\ifillnorm{U_m^i}\leq C\varepsilon\mass{U_m^i}$
 \item $\sum_{i=1}^{M_m}\mass{U_m^i} + \mass{R_m} \leq 2\mass{T_m}$
 \item $\ifillnorm{R_m}\to 0$.
\end{enumerate}
It is clear that such a decomposition yields the desired statement about $\ifillnorm{T_m}$. In order to construct such a 
decomposition one
chooses balls $B_m$ of radius $\varepsilon/2$ in such a way that each $B_m$ contains almost maximal $\|T_m\|$-measure among all balls in $X$ of radius $\varepsilon/2$.
If $\|T_m\|(B_m)$ tends to $0$ with $m\to\infty$ then it will be shown (see \propref{proposition:fillvol-lambda-sequence}) that $\ifillnorm{T_m}\to 0$. Hence the trivial 
decomposition $T_m=R_m$ satisfies the above properties.
On the other hand, if $\|T_m\|(B_m)$ stays bounded away from $0$ there are two cases to be distinguished: Either, after passing to a subsequence, all balls lie in a 
common
ball $B(y,2\varepsilon)$ or, otherwise, we may assume all balls to be disjoint. As for the first case, there exists an 
$r\in(3\varepsilon,4\varepsilon)$ such that the boundary $\bdry(T\rstr B(y,r))$ of the current obtained by restricting $T_m$ to the ball $B(y,r)$ is an integral
current and $\bdry(T\rstr B(y,r))$ is a bounded sequence converging weakly to $0$. By induction assumption there exists $S_m\in\intcurr_{k+1}(X)$ such that
$\bdry S_m=\bdry(T\rstr B(y,r))$ and $\mass{S_m}\to 0$. Using an isoperimetric inequality for cycles in $\intcurr_{k-1}(X)$ of small mass (to be proved in 
\thmref{theorem:small-mass-isoperimetric-inequality}) we show that $S_m$ can be chosen to be supported in $B(y,5\varepsilon)$. If 
$\varepsilon<\delta/5$ the $\delta$-cone type inequality yields
\begin{equation*}
 \ifillnorm{U_m^1}\leq C5\varepsilon\mass{U_m^1}
\end{equation*}
where $U_m^1:= T_m\rstr B(y,r)-S_m$.
Passing to subsequences and proceeding as above we can split off $U_m^2, U_m^3,\dots$ until ending up with a `rest' $R_m$ with $\ifillnorm{R_m}\to 0$.

\medskip

The paper is structured as follows. \secref{section:currents} contains the basic definitions and results from the theory of metric integral currents needed for the sequel. We 
point out here that our result needs no closure or compactness theorem. In \secref{section:local-cone-inequalities} we study local cone type inequalities and prove, in 
particular, that spaces for which small sets are Lipschitz contractible in a uniform way admit local cone type inequalities in all dimensions. This will in particular imply that 
Banach spaces and ${\rm CAT}(\kappa)$-spaces for all $\kappa\in\R$ admit local cone type inequalities.
The purpose of \secref{section:flat-convergence-0-dimension} is to establish \thmref{theorem:0-dimension-flat-weak}. The higher dimensional analogue, 
\thmref{theorem:flat-weak-convergence-higher-dimension}, is proved in \secref{section:flat-convergence-higher-dimension}. As mentioned above, the proof is based on 
\thmref{theorem:small-mass-isoperimetric-inequality}, also established in this section. In the last section we describe some consequences of
 \thmref{theorem:flat-weak-convergence-higher-dimension}.
%

%

%

%
%
\section{Currents in metric spaces}\label{section:currents}
The general reference for this section is \cite{Ambr-Kirch-curr} where the theory of currents in metric spaces was developed. Here, we recall those definitions
and results from \cite{Ambr-Kirch-curr} which will be needed in the sequel.

{\sloppy
Let $(X,d)$ be a complete metric space and let $\form^k(X)$ denote the set of $(k+1)$-tuples $(f,\pi_1,\dots,\pi_k)$ 
of Lipschitz functions on $X$ with $f$ bounded. The Lipschitz constant of a Lipschitz function $f$ on $X$ will
be denoted by $\lip(f)$.
}
\bd
A $k$-dimensional metric current  $T$ on $X$ is a multi-linear functional on $\form^k(X)$ satisfying the following
properties:
\begin{enumerate}
 \item If $\pi^j_i$ converges point-wise to $\pi_i$ as $j\to\infty$ and if $\sup_{i,j}\lip(\pi^j_i)<\infty$ then
       \begin{equation*}
         T(f,\pi^j_1,\dots,\pi^j_k) \longrightarrow T(f,\pi_1,\dots,\pi_k).
       \end{equation*}
 \item If $\{x\in X:f(x)\not=0\}$ is contained in the union $\bigcup_{i=1}^kB_i$ of Borel sets $B_i$ and if $\pi_i$ is constant 
       on $B_i$ then
       \begin{equation*}
         T(f,\pi_1,\dots,\pi_k)=0.
       \end{equation*}
 \item There exists a finite Borel measure $\mu$ on $X$ such that
       \begin{equation}\label{equation:mass-def}
        |T(f,\pi_1,\dots,\pi_k)|\leq \prod_{i=1}^k\lip(\pi_i)\int_X|f|d\mu
       \end{equation}
       for all $(f,\pi_1,\dots,\pi_k)\in\form^k(X)$.
\end{enumerate}
\ed
The space of $k$-dimensional metric currents on $X$ is denoted by $\curr_k(X)$ and the minimal Borel measure $\mu$
satisfying \eqref{equation:mass-def} is called mass of $T$ and written as $\|T\|$. We also call mass of $T$ the number $\|T\|(X)$ 
which we denote by $\mass{T}$.
The support of $T$ is, by definition, the closed set $\spt T$ of points $x\in X$ such that $\|T\|(B(x,r))>0$ for all $r>0$. 
Here, $B(x,r)$ denotes the closed ball $B(x,r):= \{y\in X : d(y,x)\leq r\}$.
\br
As is done in \cite{Ambr-Kirch-curr} we will also assume here that the cardinality of any set is an Ulam number. This is consistent with the 
standard ZFC set theory. We then have that $\spt T$ is separable and furthermore that $\|T\|$ is concentrated on a 
$\sigma$-compact set, i.e. $\|T\|(X\ohne C) = 0$ for a $\sigma$-compact set $C\subset X$ (see \cite{Ambr-Kirch-curr}).
\er
The restriction of $T\in\curr_k(X)$ to a Borel set $A\subset X$ is given by 
\begin{equation*}
  (T\rstr A)(f,\pi_1,\dots,\pi_k):= T(f\chi_A,\pi_1,\dots,\pi_k).
\end{equation*}
This expression is well-defined since $T$ can be extended to a functional on tuples for which the first argument lies in 
$L^\infty(X,\|T\|)$.\\
The boundary of $T\in\curr_k(X)$ is the functional
\begin{equation*}
 \bdry T(f,\pi_1,\dots,\pi_{k-1}):= T(1,f,\pi_1,\dots,\pi_{k-1}).
\end{equation*}
It is clear that $\bdry T$ satisfies conditions (i) and (ii) in the above definition. If $\bdry T$ also has 
finite mass (condition (iii)) then $T$ is called a normal current. The respective space is denoted by $\norcurr_k(X)$.\\
The push-forward of $T\in\curr_k(X)$ 
under a Lipschitz map $\varphi$ from $X$ to another complete metric space $Y$ is given by
\begin{equation*}
 \varphi_\# T(g,\tau_1,\dots,\tau_k):= T(g\circ\varphi, \tau_1\circ\varphi,\dots,\tau_k\circ\varphi)
\end{equation*}
for $(g,\tau_1,\dots,\tau_k)\in\form^k(Y)$. This defines a $k$-dimensional current on $Y$, as is easily verified.\\
In this paper we will mainly be concerned with integer rectifiable and integral currents.
For notational purposes we first repeat some well-known definitions. The Hausdorff $k$-dimensional 
measure of $A\subset X$ is defined to be
\begin{equation*}
 \hm^k(A):= \lim_{\delta\searrow 0}\inf\left\{\sum_{i=1}^\infty \omega_k\left(\frac{\diam(B_i)}{2}\right)^k :
      B\subset \bigcup_{i=1}^\infty B_i\text{, }\diam(B_i)<\delta\right\},
\end{equation*}
where $\omega_k$ denotes the Lebesgue measure of the unit ball in $\R^k$. The $k$-dimensional 
lower density $\Theta_{*k}(\mu, x)$ of a finite Borel measure $\mu$ at a point $x$ is given by the formula
\begin{equation*}
 \Theta_{*k}(\mu, x):= \liminf_{r\searrow 0}\frac{\mu(B(x,r))}{\omega_k r^k}.
\end{equation*}
An $\hm^k$-measurable set $A\subset X$
is said to be countably $\hm^k$-rectifiable if there exist countably many Lipschitz maps $f_i :B_i\longrightarrow X$ from subsets
$B_i\subset \R^k$ such that
\begin{equation*}
\hm^k(A\ohne \bigcup f_i(B_i))=0.
\end{equation*}
\bd
A current $T\in\curr_k(X)$ with $k\geq 1$ is said to be rectifiable if
\begin{enumerate}
 \item $\|T\|$ is concentrated on a countably $\hm^k$-rectifiable set and
 \item $\|T\|$ vanishes on $\hm^k$-negligible sets.
\end{enumerate}
$T$ is called integer rectifiable if, in addition, the following property holds:
\begin{enumerate}
 \item[(iii)] For any Lipschitz map $\varphi\colon X\longrightarrow \R^k$ and any open set $U\subset X$ there exists 
       $\theta\in L^1(\R^k,\Z)$ such that
       \begin{equation*}
        \varphi_\#(T\rstr U)(f,\pi_1,\dots,\pi_k)= \int_{\R^k}\theta f\det\left(\frac{\partial\pi_i}{\partial x_j}\right)d{\mathcal L}^k
       \end{equation*}
       for all $(f,\pi_1,\dots,\pi_k)\in\form^k(\R^k)$.
\end{enumerate}
\ed
A $0$-dimensional (integer) rectifiable current is a $T\in\curr_0(X)$ of the form
\begin{equation*}
 T(f)=\sum_{i=1}^\infty \theta_i f(x_i),\qquad \text{$f$ Lipschitz and bounded,}
\end{equation*}
for suitable $\theta_i\in\R$ (or $\theta_i\in \Z$) and $x_i\in X$.\\
The space of rectifiable currents is denoted by $\rectcurr_k(X)$, that of integer rectifiable currents by $\intrectcurr_k(X)$.
Endowed with the mass norm $\curr_k(X)$ is a Banach space, $\rectcurr_k(X)$ a closed subspace, and $\intrectcurr_k(X)$ a closed additive subgroup. 
This follows directly from the definitions.
Integer rectifiable normal currents are called integral currents. The respective space is denoted by $\intcurr_k(X)$.
As the mass of a $k$-dimensional normal current vanishes on $\hm^k$-negligible sets (\cite[Theorem 3.7]{Ambr-Kirch-curr}) it 
is easily verified that the push-forward of an integral current under a Lipschitz map is again an integral current. 
In the following, an element $T\in\intcurr_k(X)$ with zero boundary $\bdry T=0$ will be called a cycle. An element $S\in\intcurr_{k+1}(X)$
satisfying $\bdry S=T$ is said to be a filling of $T$.\\
The characteristic set $S_T$ of a rectifiable current $T\in\rectcurr_k(X)$ is defined by
\begin{equation}\label{equation:minimal-set}
 S_T:= \{x\in X: \Theta_{\star k}(\|T\|, x)>0\}.
\end{equation}
It can be shown that $S_T$ is countably $\hm^k$-rectifiable and that $\|T\|$ is concentrated on $S_T$.
In the next theorem the function $\lambda: S_T\longrightarrow (0,\infty)$ denotes the area factor on the (weak) tangent spaces to $S_T$
as defined in \cite{Ambr-Kirch-curr}. We do not provide a definition here since for our purposes it is enough to know 
that $\lambda$ is $\hm^k$-integrable and bounded from below by $k^{-k/2}$ (see \cite[Lemma 9.2]{Ambr-Kirch-curr}).
\begin{thm}[{\cite[Theorem 9.5]{Ambr-Kirch-curr}}]\label{theorem:mass-representation}
 If $T\in\rectcurr_k(X)$ then there exists a $\hm^k$-integrable function $\theta: S_T\longrightarrow (0,\infty)$ such that
 \begin{equation*}
  \|T\|(A)=\int_{A\cap S_T}\lambda\theta d\hm^k\qquad\text{for $A\subset X$ Borel,}
 \end{equation*}
 that is, $\|T\|= \lambda\theta d\hm^k\rstr S_T$.
 Moreover, if $T$ is an integral current then $\theta$ takes values in $\N:=\{1,2,\dots\}$ only.
\end{thm}
\thmref{theorem:mass-representation} yields a lower density bound around almost every point in the characteristic set of an arbitrary non-trivial integral current. 
We will use this in the proofs of \thmref{theorem:small-mass-isoperimetric-inequality} and \propref{proposition:reduction-density}.

The main technique we will use in the proof of Theorems \ref{theorem:0-dimension-flat-weak} and \ref{theorem:flat-weak-convergence-higher-dimension} is that 
of slicing. We will, in particular, make frequent use of the following important result known as the Slicing Theorem and proved in 
\cite[Theorems 5.6 and 5.7]{Ambr-Kirch-curr}.
\bt\label{theorem:slicing}
Let be $T\in\norcurr_k(X)$ and $\varrho$ a Lipschitz function on $X$. Then there exists for almost every $r\in\R$ a normal current
$\slice{T}{\varrho}{r}\in\norcurr_{k-1}(X)$ with the following properties:
\begin{enumerate}
 \item $\slice{T}{\varrho}{r}= \bdry(T\rstr\{\varrho \leq r\}) - (\bdry T)\rstr\{\varrho\leq r\}$
 \item $\|\slice{T}{\varrho}{r}\|$ and $\|\bdry\slice{T}{\varrho}{r}\|$ are concentrated on $\varrho^{-1}(\{r\})$
 \item $\mass{\slice{T}{\varrho}{r}}\leq\lip(\varrho)\frac{d}{dr}\mass{T\rstr\{\varrho\leq r\}}$.
\end{enumerate}
Moreover, if $T\in\intcurr_k(X)$ then $\slice{T}{\varrho}{r}\in\intcurr_{k-1}(X)$ for almost all $r\in\R$.
\et
%

\section{Spaces admitting local cone type inequalities}\label{section:local-cone-inequalities}

In this section we establish local and global cone type inequalities for various classes of spaces.
In the first paragraph we define the product of a current with an interval. In the subsequent paragraph we show that spaces in which balls are 
null-homotopic in a Lipschitz way (the Lipschitz constant depending on the diameter) admit a cone type inequality. 
From this we immediately obtain cone type inequalities for Banach and ${\rm CAT}(\kappa)$-spaces for all $\kappa\in\R$.
\subsection{Products of currents}
The following construction is a slightly modified version of the one given in the first part of Section 10 of \cite{Ambr-Kirch-curr}.\\
Let $(X,d)$ be a complete metric space and endow $[0,1]\times X$ with the Euclidean product metric.
Given a Lipschitz function $f$ on $[0,1]\times X$ and $t\in[0,1]$ we define the function $f_t:X\longrightarrow \R$ by
$f_t(x):= f(t,x)$. To every $T\in\norcurr_k(X)$, $k\geq 1$, and every $t\in[0,1]$ we associate the normal $k$-current on $[0,1]\times X$ given
by the formula
\begin{equation*}
  ([t]\times T)(f,\pi_1,\dots,\pi_k):= T(f_{t},\pi_{1\,t},\dots,\pi_{k\,t}).
\end{equation*}
The product of a normal current with the interval $[0,1]$ is defined 
as follows.
\sloppy
 \bd
  For a normal current $T\in\norcurr_k(X)$
  the functional $[0,1]\times T$ on $\form^{k+1}([0,1]\times X)$ is given by
  \begin{equation*}
  \begin{split}
   ([0,1]\times& T) (f,\pi_1,\dots,\pi_{k+1}):= \\
      &\sum_{i=1}^{k+1}(-1)^{i+1}\int_0^1T\left(f_t\frac{\partial \pi_{i\,t}}{\partial t},\pi_{1\,t},
                                               \dots,\pi_{i-1\,t},\pi_{i+1\,t},\dots,\pi_{k+1\,t}\right)dt
  \end{split}
  \end{equation*}
  for $(f,\pi_1,\dots,\pi_{k+1})\in\form^{k+1}([0,1]\times X)$.
 \ed
\fussy
We have the following result whose proof is analogous to that
of \cite[Proposition 10.2 and Theorem 10.4]{Ambr-Kirch-curr}.
\bt\label{theorem:cone-construction}
  For every $T\in\norcurr_k(X)$, $k\geq 1$, with bounded support the functional $[0,1]\times T$ is a $(k+1)$-dimensional normal current
 on $[0,1]\times X$ with boundary
  \begin{equation*}
   \partial([0,1]\times T)= [1]\times T - [0]\times T - [0,1]\times\partial T.
  \end{equation*}
 Moreover, if $T\in\intcurr_k(X)$ then $[0,1]\times T\in\intcurr_{k+1}([0,1]\times X)$.
\et
%

\subsection{Cone type inequalities and $\gamma$-Lipschitz contractibility}\label{subsection:cone-lipschitz-contractibility}
Let $(X,d)$ be a metric space. We say that a bounded subset $B\subset X$ is $(\beta,\gamma)$-Lipschitz contractible in $X$ if there exists a map  
 $\varphi: [0,1]\times B\to X$ satisfying $\varphi(1,\cdot)=\operatorname{id}_B$ and $\varphi(0,\cdot)\equiv x_0$ for some $x_0\in X$
 and moreover
 \begin{equation*}
  d(\varphi(t,x),\varphi(t',x'))\leq \beta |t-t'| + \gamma d(x,x')\qquad\text{for all $x,x'\in B$ and $t,t'\in[0,1]$}.
 \end{equation*}
A bounded subset $B\subset X$ which is $(\gamma\diam B,\gamma)$-Lipschitz contractible will be called $\gamma$-Lipschitz contractible and the
corresponding map $\varphi$ a $\gamma$-contraction of $B$. 
%
%
\br
 The property that every bounded subset is $\gamma$-contractible is invariant under bi-Lipschitz homeomorphisms in the following sense: 
 If $\varphi: X\to Y$ is a bi-Lipschitz homeomorphism then bounded subsets in $X$ are $\gamma$-Lipschitz contractible for some $\gamma$ if and only if bounded
 subsets in $Y$ are $\gamma'$-Lipschitz contractible for some $\gamma'$.
\er
We show that spaces all of whose subsets of small diameter are $\gamma$-Lipschitz contractible for a fixed $\gamma$ admit local cone type inequalities for $\intcurr_k(X)$ 
for every $k\geq 1$.
The `cone filling' of $T\in\intcurr_k(X)$ is constructed as follows: 
Let $\varphi:[0,1]\times\spt T\to X$ be a $\gamma$-contraction of $\spt T$. 
Since $\spt([0,1]\times T)=[0,1]\times\spt T$ and since $\varphi$ is Lipschitz
also with respect to the Euclidean product metric on $[0,1]\times\spt T$ it follows that $S:= \varphi_\#([0,1]\times T)$ is well-defined and, 
by \thmref{theorem:cone-construction}, $S\in\intcurr_{k+1}(X)$ with $\bdry S=T$ if $\bdry T=0$.
%
 \bp\label{proposition:coneineq-locally-contractible}
  If $(X,d)$ is a complete metric space all of whose subsets of diameter no larger than $\delta$ are $\gamma$-Lipschitz contractible then for every 
  cycle $T\in\intcurr_k(X)$, $k\geq 1$, satisfying $\diam(\spt T)\leq\delta$ there exists an $S\in\intcurr_{k+1}(X)$ satisfying $\bdry S=T$ and
  \begin{equation*}
   \mass{S}\leq (k+1)\gamma^{k+1}\diam(\spt T)\mass{T},
  \end{equation*}
  i.e. $X$ admits a $\delta$-cone type inequality for $\intcurr_k(X)$ with $k\geq 1$.
 \ep
We remark that a space $X$ satisfying the assumptions of the proposition is automatically $\delta$-quasiconvex.
 \begin{proof}
  Let $T$, $\varphi$, and $S$ be as above. Then for fixed $x\in\spt T$ the map $t\mapsto \varphi(t,x)$ is 
  $\gamma\diam(\spt T)$-Lipschitz, whereas for fixed $t\in[0,1]$ 
  the map $x\mapsto \varphi(t,x)$ is $\gamma$-Lipschitz. 
 For $(f,\pi_1,\dots,\pi_{k+1})\in\form^{k+1}(X)$ we therefore obtain
  \begin{equation*}
   \begin{split}
    |&S(f,\pi_1,\dots,\pi_{k+1})| \\
       &\leq \sum_{i=1}^{k+1}\left|\int_0^1 
             T\left(f\circ\varphi_t\frac{\partial (\pi_i\circ\varphi_t)}{\partial t},\pi_1\circ\varphi_t,\dots,\pi_{i-1}\circ\varphi_t,
               \pi_{i+1}\circ\varphi_t,\dots,\pi_{k+1}\circ\varphi_t\right)dt\right|\\
                                        &\leq \sum_{i=1}^{k+1} \int_0^1\prod_{j\not=i}\lip(\pi_j\circ\varphi_t)
                                                \int_X\left|f\circ\varphi_t\;\frac{\partial (\pi_i\circ\varphi_t)}{\partial t}\right|
                                                 d\|T\|dt\\
                                        &\leq (k+1)\gamma^{k+1}\diam(\spt T)\prod_{j=1}^{k+1}\lip(\pi_j)\int_0^1
                                                \int_X|f\circ\varphi(t,x)|d\|T\|(x)dt.
   \end{split}
  \end{equation*}
  From this it follows that $\|S\|\leq (k+1)\gamma^{k+1}\diam(\spt T)\varphi_{\#}({\mathcal L}^1\times\|T\|)$ and 
  this concludes the proof.
 \end{proof} 
An important class of examples satisfying the assumptions of the proposition above is given by Banach spaces and ${\rm CAT}(\kappa)$-spaces. These are complete 
geodesic metric spaces for which all geodesic triangles of perimeter strictly smaller than $2\pi/\sqrt{\kappa}$ are at least as slim as corresponding triangles
in the model space of curvature $\kappa$. This condition is some kind of global curvature condition. For example, complete Riemannian manifolds of sectional curvature
bounded by $\kappa$ are locally ${\rm CAT}(\kappa)$. 
See e.g.~\cite{Ballmann}, \cite{Bridson-Haefliger}, \cite{Burago-Burago-Ivanov} for an account on such spaces. If $X$ is a Banach space
and $B\subset X$ a bounded subset then we can simply define $\varphi(t,x):= tx+(1-t)x_0$ for $x\in B$ where $x_0\in B$ is an 
arbitrary given point. This clearly defines a $1$-contraction of $B$ by the triangle inequality. 
Analogously, if $X$ is a ${\rm CAT}$-space then for sets $B\subset X$ with small diameter (depending on the curvature bound on $X$) one sets $\varphi(t,x):= c_{x_0x}(t)$ where $c_{x_0x}$ is the 
unique reparameterized geodesic from $x_0$ to $x$. Triangle comparison then yields immediately that this defines a Lipschitz contraction. Further examples of spaces satisfying the assumptions
of \propref{proposition:coneineq-locally-contractible} are given by spaces of non-positive curvature in the sense of Busemann, compact neighborhood Lipschitz retracts in $\R^n$, Lipschitz retracts
of spaces admitting local cone type inequalities, etc.

\section{Flat convergence for $0$-dimensional integral currents}\label{section:flat-convergence-0-dimension}

Let $(X,d)$ be $\delta$-quasiconvex and $T\in\intcurr_0(X)$ such that $\diam(\spt T)\leq\delta$. 
Then there exists an $S\in\intcurr_1(X)$ with $\bdry S=T$ if and only if 
$T(1)=0$. In this case an $S$ which will suit our purposes can be constructed as follows. 
Write $T$ as $T=\sum_{i=1}^{M}\Lbrack x_i^+\Rbrack - \Lbrack x_i^-\Rbrack$ in such a way that $2M=\mass{T}$. Here we have set
$\Lbrack x\Rbrack(f):=f(x)$ for $f\in\lip_b(X)$. Then we define
\begin{equation*}
 S:= \sum_{i=1}^{M} (\gamma_{x_i^-x_i^+})_{\#}\Lbrack \chi_{[0,1]}\Rbrack,
\end{equation*}
where $\gamma_{x_i^-x_i^+}$ is as in \defref{definition:quasi-convex} and parameterized in such a way that it has Lipschitz constant at most $Cd(x_i^-,x_i^+)$. 
Then, clearly, $\bdry S=T$ and moreover
\begin{equation*}
 \mass{S}\leq C\diam(\spt T)M=\frac{1}{2}C\diam(\spt T)\mass{T}.
\end{equation*}
In the proof of Theorems \ref{theorem:0-dimension-flat-weak} and \ref{theorem:flat-weak-convergence-higher-dimension}
we will need the following auxiliary result which is a reformulation of \cite[Proposition 8.3]{Ambr-Kirch-curr}.
\bl\label{lemma:boundary-slices-sequence}
 Let $k\geq 0$ and let $(T_m)\subset\intcurr_k(X)$ be a bounded sequence of cycles weakly converging to $0$. Let furthermore $A\subset X$ be an arbitrary subset and denote 
 by $A(r)$ the closed $r$-neighborhood of $A$ whenever $r>0$. Then for almost every $r>0$ there exists a subsequence $(T_{m_j})$ such that
 \begin{enumerate}
  \item $T_{m_j}\rstr A(r)\in\intcurr_k(X)$ for every $j\in\N$
  \item $T_{m_j}\rstr A(r)$ and hence also $\bdry(T_{m_j}\rstr A(r))$ weakly converges to $0$ as $j\to\infty$
  \item $\mass{\bdry(T_{m_j}\rstr A(r))}\leq C$ for every $j\in\N$ for a constant $C=C(r)<\infty$.
 \end{enumerate}
\el
In case $k=0$, of course $\bdry(T_{m_j}\rstr A(r))$ is not defined and the statements in (ii) and (iii) concerning $\bdry(T_{m_j}\rstr A(r))$ should be
omitted.

We are now in a position to prove the equivalence of flat and weak convergence for $0$-dimensional integral currents.

\begin{proof}[{Proof of \thmref{theorem:0-dimension-flat-weak}}]
 It is clearly enough to consider the case $T=0$ and to find a subsequence $T_{m_j}$ for which $\ifillnorm{T_{m_j}}\to 0$ as $j\to\infty$. 
 After passing to a subsequence we may assume the existence of $M\in\N$, $a_1,\dots,a_M\in\Z\ohne\{0\}$, and $x_m^1,\dots,x_m^M\in X$ such that 
 \begin{equation*}
  T_m(f) = \sum_{i=1}^M a_i f(x_m^i)
 \end{equation*}
 for all $f\in\lip_b(X)$ and all $m\in\N$.
 Let $0<\varepsilon\leq\delta/6$ be arbitrary. We claim that after passing to a subsequence there exists for every $m\in\N$ a decomposition
 $T_m=T'_m+R_m$ for some $T'_m,R_m\in\intcurr_0(X)$ such that 
 \begin{enumerate}
  \item $\mass{T'_m}\geq 1$
  \item $\mass{T_m}=\mass{T'_m}+\mass{R_m}$
  \item $\ifillnorm{T'_m}\leq 6C\varepsilon\mass{T'_m}$
 \end{enumerate}
 and
 \begin{enumerate}
 \setcounter{enumi}{3}
  \item $R_m\to 0$ as $m\to\infty$.
 \end{enumerate}
 The lemma clearly follows by successive application of this claim and by letting $\varepsilon\to 0$.
 In order to prove the claim we distinguish two cases. First, suppose that after passage to a subsequence there exist indices $i_m\in\{1,\dots,M\}$ such that 
 $x_m^{i_m}\in B(x_1^{i_1},2\varepsilon)$ for all $m\in\N$. By \lemref{lemma:boundary-slices-sequence} there exists $r\in(2\varepsilon, 3\varepsilon]$ such 
 that after further 
 passage to a subsequence (also for the $i_m$) we obtain $T'_m:= T_m\rstr B(x_1^{i_1},r)\to 0$ and thus also $R_m:=T_m\rstr (X\ohne B(x_1^{i_1},r))\to 0$ as 
 $m\to \infty$, and hence (iv). 
 It is clear that $T'_m,R_m\in\intcurr_0(X)$ and $T_m=T'_m+R_m$ and that furthermore assertion (ii) holds. 
 Since $T'_m\to 0$ we have in particular that $T'_m(1)\to0$ and 
 therefore $T'_m(1)=0$ for $m$ large enough.
 Since $\diam(\spt T'_m)\leq2r\leq6\varepsilon\leq\delta$ we can use the construction at the beginning of this section to obtain 
 $\ifillnorm{T'_m}\leq 3C\varepsilon\mass{T'_m}$ which proves (iii). Property (i) is clear since 
 $d(x_m^{i_m},x_1^{i_1})\leq 2\varepsilon$ and hence $\mass{T'_m}\geq |a_{i_m}|\geq 1$. We now consider the second case for which we 
 assume that no such subsequence and indices as in the first case exist. Then for every $m\in\N$ and every $i\in\{1,\dots,M\}$ there exist only finitely 
 many $m'$ and $i'$ with 
 $x_{m'}^{i'}\in B(x_m^i,2\varepsilon)$. After passing to a subsequence we may therefore assume the sets $E_m(\varepsilon):= \cup_{i=1}^M B(x_m^i,\varepsilon)$ 
 to be pairwise disjoint. Define $A:= \{x_m^1: m\in\N\}$. By construction we have
 $T_m\rstr A(\varepsilon')=T_m\rstr B(x_m^1,\varepsilon')$ whenever $0<\varepsilon'\leq\varepsilon$. By \lemref{lemma:boundary-slices-sequence} there exists
 $0<\varepsilon'<\varepsilon$ arbitrarily close to $\varepsilon$ such that after passage to a subsequence we have $T'_m:= T_m\rstr B(x_m^1,\varepsilon')\to 0$ 
 and hence also $R_m:= T_m-T'_m\to 0$ as $m\to\infty$. 
 Assertions (i) to (iv) then follow as in the first case.
 This proves the claim above and therefore the lemma.
\end{proof}

\section{Flat convergence for $k$-dimensional integral currents for $k\geq 1$}\label{section:flat-convergence-higher-dimension}

\subsection{Splitting off cycles of controlled diameter}
We start with a simple but useful observation.
 \bl\label{lemma:polynomial}
  Fix $\bar{C}>0$, $k\geq 2$, $0\leq r_0<r_1<\infty$, and suppose $\beta : [r_0,r_1]\longrightarrow (0,\infty)$ is non-decreasing and satisfies
  \begin{enumerate}
   \item $\beta(r_0)= \frac{r_0^k}{\bar{C}^{k-1}k^k}$
   \item $\beta(r)\leq \bar{C}[\beta'(r)]^{k/(k-1)}$ for almost every $r\in(r_0,r_1)$.
  \end{enumerate}
  Then it follows that
  \begin{equation*}
   \beta(r)\geq \frac{r^k}{\bar{C}^{k-1}k^k}\quad\text{ for all } r\in[r_0,r_1].
  \end{equation*}
 \el
\begin{proof}
 By rearranging (ii) we obtain
 \begin{equation*}
  \frac{\beta'(t)}{\beta(t)^{\frac{k-1}{k}}}\geq \frac{1}{\bar{C}^{\frac{k-1}{k}}}
 \end{equation*}
 and integration from $r_0$ to $r$ yields the claimed estimate.
\end{proof}
\bd
 Let $(X,d)$ be a complete metric space and $k\in\N$. Then $X$ is said to admit an isoperimetric inequality of Euclidean type 
 with constant $C$ for cycles in $\intcurr_k(X)$ with mass no larger than $C'$ if for every cycle $T\in\intcurr_k(X)$ with $\mass{T}\leq C'$ there 
 exists an $S\in\intcurr_{k+1}(X)$ with $\bdry S=T$ and such that
\begin{equation*}
 \mass{S}\leq C[\mass{T}]^{\frac{k+1}{k}}.
\end{equation*}
\ed
The next result roughly states that cycles of small mass in spaces admitting an isoperimetric inequality of Euclidean type for cycles of small mass 
have fillings which stay near their boundary.
\bl\label{lemma:local-isoperimetric-support}
 Let $(X,d)$ be a complete metric space and $k\geq 2$. Suppose $X$ admits an isoperimetric inequality of Euclidean type 
 with constant $C$ for cycles in $\intcurr_{k-1}(X)$ of mass at most $C'$. Then for every $T\in\intcurr_{k-1}(X)$ with $\bdry T=0$ and $\mass{T}\leq C'$ 
 there exists an $S\in\intcurr_k(X)$ satisfying $\bdry S=T$ as well as
 \begin{equation*}
  \mass{S}\leq C[\mass{T}]^{\frac{k}{k-1}}
 \end{equation*}
 and
 \begin{equation*}
  \spt S\subset B\left(\spt T, 3Ck\mass{T}^{\frac{1}{k-1}}\right).
 \end{equation*}
\el
The proof is a variation of the argument in \cite[Theorem 10.6]{Ambr-Kirch-curr}.
\begin{proof}
Let ${\mathcal M}$ denote the space consisting of all $S\in\intcurr_k(X)$ with $\bdry S=T$ and endow ${\mathcal M}$ with the metric
$d_{\mathcal M}(S,S'):=\mass{S-S'}$. Then ${\mathcal M}$ is a complete metric space. Choose an $\tilde{S}\in{\mathcal M}$ satisfying
$\mass{\tilde{S}}\leq C[\mass{T}]^{k/(k-1)}$. By the Ekeland variation principle \cite{Ekeland}  there exists an $S\in{\mathcal M}$ with $\mass{S}\leq \mass{\tilde{S}}$ 
and such that the function 
\begin{equation*}
 S'\mapsto \mass{S'}+\frac{1}{2}\mass{S-S'}
\end{equation*}
is minimal at $S'=S$. Let be $x\in\spt S\ohne\spt T$ and set $\varrho_x(y):= d(x,y)$ and $\beta(r):= \|S\|(B(x,r))$. 
Then, \thmref{theorem:slicing} implies that for almost every $0<r<\dist(x,\spt T)$ the slice $\slice{S}{\varrho_x}{r}$ exists, has zero boundary, 
belongs to $\intcurr_{k-1}(X)$ and furthermore satisfies $\mass{\slice{S}{\varrho_x}{r}}\leq \beta'(r)$. 
We claim that 
\begin{equation*}
 \beta(r)\geq \frac{r^k}{(3C)^{k-1}k^k}\quad\text{for all $0\leq r \leq\dist(x,\spt T)$}.
\end{equation*}
Indeed, if $\mass{\slice{S}{\varrho_x}{r}}\leq C'$ then we choose an
$S_r\in\intcurr_k(X)$ with $\bdry S_r=\slice{S}{\varrho_x}{r}$ and $\mass{S_r}\leq C[\mass{\slice{S}{\varrho_x}{r}}]^{k/(k-1)}$ which exists by assumption. 
The integral current
$S\rstr (X\ohne B(x,r))+S_r$ has boundary $T$ and thus, comparison with $S$ yields
\begin{equation*}
\mass{S\rstr (X\ohne B(x,r))+S_r} + \frac{1}{2}\mass{S\rstr B(x,r)-S_r}\geq \mass{S}.
\end{equation*}
This together with the isoperimetric inequality yields
\begin{equation*}
 \beta(r)=\mass{S\rstr B(x,r)}\leq3\mass{S_r}\leq 3C[\mass{\slice{S}{\varrho_x}{r}}]^{\frac{k}{k-1}}\leq 3C[\beta'(r)]^{\frac{k}{k-1}}.
\end{equation*}
If, on the other hand, $\mass{\slice{S}{\varrho_x}{r}}> C'$ then
\begin{equation*}
 \beta(r)\leq \mass{S}\leq C[\mass{T}]^{\frac{k}{k-1}}\leq C(C')^{\frac{k}{k-1}}<C[\mass{\slice{S}{\varrho_x}{r}}]^{\frac{k}{k-1}}\leq C[\beta'(r)]^{\frac{k}{k-1}},
\end{equation*}
and thus $\beta(r)\leq 3C[\beta'(r)]^{k/(k-1)}$ for almost every $0<r<\dist(x,\spt T)$.
The claim now follows from \lemref{lemma:polynomial}. The second estimate in the lemma is a direct consequence of this claim.
\end{proof}
The proof of \lemref{lemma:local-isoperimetric-support} can also be used to show that if $k\geq 2$, if $(X,d)$ is as in \lemref{lemma:local-isoperimetric-support}, 
and if $T\in\intcurr_{k-1}(X)$ is a cycle with
$\ifillnorm{T}<C(C')^{k/(k-1)}$ (no restriction on the mass of $T$), then
\begin{equation*}
 \fillrad{T}\leq (3C)^{\frac{k-1}{k}}k\ifillnorm{T}^{\frac{1}{k}}.
\end{equation*}
Here, $\fillrad{T}$ is defined by
\begin{equation*}
 \fillrad{T}:= \inf\{r\geq0: \exists S'\in\intcurr_k(X),\, \bdry S'=T,\; \spt S'\subset B(\spt T,r)\}.
\end{equation*}
\bl\label{lemma:partial-cover}
  Let $(Y,d)$ be a metric space, $\mu$ a finite Borel measure on $Y$, and $F>0$, $k\in\N$. Suppose $r_0:Y\to[0,\infty)$ is a function such that for 
  $\mu$-almost every $y\in Y$ we have $r_0(y)>0$ and
  \begin{enumerate}
   \item $\mu(B(y,r_0(y)))= F[r_0(y)]^k$
   \item $\mu(B(y,5r_0(y)))<5^kF[r_0(y)]^k$.
  \end{enumerate}
  Then there exist finitely many points $y_1,\dots, y_N\in Y$ satisfying
  \begin{enumerate}
   \item $r_0(y_i)>0$
   \item $B\left(y_i, 2r_0(y_i)\right)\cap B\left(y_j, 2r_0(y_j)\right)=\emptyset$ if $i\not= j$
   \item $\sum_{i=1}^N\mu(B(y_i,r_0(y_i)))\geq \frac{1}{5^k}\mu(Y)$.
  \end{enumerate}
\el
This statement is analogous to that of Lemma 3.2 in \cite{Wenger-GAFA}, the proof of which is almost that of the simple Vitali lemma.
\bl\label{lemma:general-decomposition}
 Let $(X,d)$ be a complete metric space and $k\in\N$. If $k\geq 2$ then suppose furthermore that $X$ admits an
 isoperimetric inequality of Euclidean type with a constant $C>0$ for cycles in $\intcurr_{k-1}(X)$ of mass no larger than $C'$.
 Let $\lambda\in(0,1/6)$ and $T\in\intcurr_k(X)$ with $\bdry T=0$ and suppose $r_0,r_1:X\to[0,\infty)$ are functions such that for $\|T\|$-almost every 
 $y\in X$ we have
 \begin{enumerate}
  \item[(a)] $0<r_0(y)<r_1(y)\leq \frac{4}{3}r_0(y)$\\*[-1ex]
  \item[(b)] $\|T\|(B(y,r_0(y)))= F[r_0(y)]^k$\\*[-1ex]
  \item[(c)] $\|T\|(B(y,r_1(y)))<F[r_1(y)]^k$\\*[-1ex]
  \item[(d)] $\|T\|(B(y,5r_0(y)))<\min\{G, 5^kF[r_0(y)]^k\}$,
 \end{enumerate}
 where $F:=1$  and $G:=\infty$ if $k=1$ and $F:= \lambda^{k-1}C^{1-k}k^{-k}$ and $G:=C(C')^\frac{k}{k-1}\lambda^{-1}$ if $k\geq 2$.
 Then there exist finitely many points $y_1,\dots,y_N\in X$ and a decomposition $$T=T_1+\dots+T_N+R$$ such that $T_i,R\in\intcurr_k(X)$ are cycles with
 the following properties:
 \begin{enumerate}
  \item $\spt T_i\subset B(y_i,2r_0(y_i))$ and 
   \begin{equation*}
    \diam(\spt T_i)\leq \frac{4}{[F(1-\lambda)]^{\frac{1}{k}}}[\mass{T_i}]^{\frac{1}{k}}
   \end{equation*}
  \item $\sum_{i=1}^N\mass{T_i}\leq (1+\lambda)\mass{T}$\\*[-1ex]
  \item $\mass{R}\leq \left(1-\frac{1}{5^k}(1-\lambda)\right)\mass{T}$\\*[-1ex]
  \item $\|R\|(B(y,\varepsilon))\leq \|T\|(B(y,\varepsilon))+\lambda\mass{T}$ for all $y\in X$ and all $\varepsilon>0$.
 \end{enumerate}
\el
\begin{proof}
 Clearly, $\mu:=\|T\|$ and $r_0$ satisfy the conditions of \lemref{lemma:partial-cover} and therefore there exist points
 $y_1,\dots,y_N\in X$ such that $r_0(y_i)>0$ and such that the $B(y_i,2r_0(y_i))$ are pairwise disjoint and 
 \begin{equation}\label{equation:covering-sum-lower-bound}
  \sum_{i=1}^N\|T\|(B(y_i,r_0(y_i)))\geq \frac{1}{5^k}\mass{T}.
 \end{equation}
 Fix $i\in\{1,\dots,N\}$ and abbreviate (for the moment) $r_0:= r_0(y_{i})$ as well as $r_1:=r_1(y_{i})$. Define $\beta(r):= \|T\|(B(y_{i},r))$ and
 observe that $\beta$ is non-decreasing and satisfies $\beta(r_0)=Fr_0^k$ and $\beta(r_1)<Fr_1^k$.
 Denoting by $\varrho$ the function $\varrho(x):= d(y_i,x)$, \thmref{theorem:slicing} implies that the slice 
 $\slice{T}{\varrho}{r}=\bdry(T\rstr B(y_i,r))$ exists for almost all $r$, is a cycle in $\intcurr_{k-1}(X)$, and satisfies moreover
 \begin{equation}\label{equation:slice-derivative}
  \mass{\slice{T}{\varrho}{r}}\leq \beta'(r)\quad\text{for almost every $r>0$}.
 \end{equation}
 We now consider one dimensional and higher dimensional cycles separately: If $k=1$ then $F=1$ and thus, by (b) and (c), there exists a measurable set 
 $\Omega\subset[r_0,r_1]$ of positive measure such that $\beta'(r)<1$ for all $r\in\Omega$. By \thmref{theorem:slicing}, we may assume that for every $r\in\Omega$ 
 the slice $\slice{T}{\varrho}{r}$ exists and is a 
 $0$-dimensional integral current. Therefore, $\mass{\slice{T}{\varrho}{r}}$ is an integer number and hence it follows from 
 \eqref{equation:slice-derivative} that $\mass{\slice{T}{\varrho}{r}}=0$. Therefore $T_{i}:=T\rstr B(y_{i},r)$ satisfies $\bdry T_{i}=0$ and by (a)
 \begin{equation*}
  \diam(\spt T_{i})\leq 4r_0\left(y_{i}\right)=4\|T\|(B(y_{i},r_0(y_{i})))\leq 4\mass{T_{i}}.
 \end{equation*}
 This proves assertion (i) in the case $k=1$.

 If $k\geq 2$ then \lemref{lemma:polynomial} and conditions (b) and (c) imply the existence of $\Omega\subset[r_0,r_1]$ 
 of positive measure such that
\begin{equation}\label{eq:beta-growth}
 C[\beta'(r)]^{\frac{k}{k-1}}<\lambda\beta(r)\quad\text{for all $r\in\Omega$}.
\end{equation}
By \thmref{theorem:slicing} we may assume without loss of generality that the slice $\slice{T}{\varrho}{r}$ exists for every
$r\in\Omega$, is a cycle in $\intcurr_{k-1}(X)$ and satisfies $\mass{\slice{T}{\varrho}{r}}\leq \beta'(r)$. Using (d) we therefore obtain 
\begin{equation*}
 \mass{\slice{T}{\varrho}{r}}\leq \beta'(r)\leq\left(\frac{\lambda}{C}\right)^{\frac{k-1}{k}}[\beta(r)]^{\frac{k-1}{k}}< C'
\end{equation*}
for every $r\in\Omega$
and hence, by assumption and by \lemref{lemma:local-isoperimetric-support}, there exists an $S\in\intcurr_k(X)$ with $\bdry S=\slice{T}{\varrho}{r}$ 
which satisfies 
\begin{equation*}
 \mass{S}\leq C[\mass{\slice{T}{\varrho}{r}}]^{\frac{k}{k-1}}
\end{equation*}
and moreover 
\begin{equation}\label{equation:support-in-ball}
 \spt S\subset B\left(\spt\slice{T}{\varrho}{r},3Ck\mass{\slice{T}{\varrho}{r}}^{\frac{1}{k-1}}\right).
\end{equation}
In particular, from \eqref{equation:slice-derivative} and \eqref{eq:beta-growth}, we conclude
\begin{equation}\label{eq:fill-estimate}
\mass{S}\leq\lambda\beta(r).
\end{equation}
It follows that $T_i:=T\rstr B(y_i,r) - S\in\intcurr_k(X)$ has zero boundary and satisfies
\begin{equation}\label{equation:lower-upper-mass-bound-cycles}
 (1-\lambda)\|T\|(B(y_i,r))\leq \mass{T_i}\leq(1+\lambda)\|T\|(B(y_i,r)).
\end{equation}
Moreover, \eqref{equation:slice-derivative}, \eqref{eq:beta-growth}, and \eqref{equation:support-in-ball} yield
\begin{equation*}
 \spt S\subset B\left(y_i, r + 3C^{\frac{k-1}{k}}\lambda^{\frac{1}{k}}k[\|T\|(B(y_i,r))]^{\frac{1}{k}}\right)
\end{equation*}
and hence, using the definition of $F$, condition (c) and 
the fact that $\lambda\leq \frac{1}{6}$ and $r\leq \frac{4}{3}r_0(y_i)$, we obtain that 
$T_i$ has support in a ball centered at $y_i$ with radius $\bar{r}$ satisfying
\begin{equation*}
 \bar{r}\leq r+3C^{\frac{k-1}{k}}\lambda^{\frac{1}{k}}k[\|T\|(B(y_i,r))]^{\frac{1}{k}}\leq 2r_0(y_i).
\end{equation*}
From this we conclude
\begin{equation}\label{equation:support-decomposition}
\diam(\spt T_i)\leq 4r_0(y_i) =  \frac{4}{F^{\frac{1}{k}}}[\beta(r_0(y_i))]^\frac{1}{k} 
               \leq \frac{4}{[F(1-\lambda)]^\frac{1}{k}}\mass{T_i}^\frac{1}{k}
\end{equation}
which proves assertion (i) in the case $k\geq 2$. 

Since our construction of $T_i$ leaves $T\rstr (X\ohne B(y_i,2r_0(y_i)))$ unaffected (by the fact that the balls
$B(y_i,2r_0(y_i))$ are pairwise disjoint) we can apply the above construction for every $i\in\{1,\dots,N\}$ to obtain cycles $T_1,\dots,T_N$.
Setting $R:=T-\sum_{i=1}^NT_i$ we obtain a decomposition satisfying the claimed properties. 
Indeed, by \eqref{equation:lower-upper-mass-bound-cycles},
\begin{equation*}
 \sum_{i=1}^N\mass{T_i}\leq (1+\lambda)\sum_{i=1}^N \|T\|(B(y_i,r_i)\leq (1+\lambda)\mass{T}
\end{equation*}
where $r_i$ is the particular $r$ chosen in the above decomposition procedure for the index $i$. This proves assertion (ii). 
Using \eqref{equation:covering-sum-lower-bound} we compute 
\begin{equation}\label{equation:mass-rest-cycle}
 \mass{R}\leq \|T\|\left(X\ohne \bigcup B(y_i,r_i)\right) + \lambda \sum\|T\|(B(y_i,r_i))
         \leq (1-\alpha(1-\lambda))\mass{T}
\end{equation}
and furthermore
\begin{equation*}
 \begin{split}
  \|R\|(B(y,\varepsilon))&\leq \|T\|\left(B(y,\epsilon)\ohne\bigcup B(y_i,r_i)\right)+\sum_{i=1}^N\|S_i\|(B(y,\varepsilon))\\
               &\leq \|T\|(B(y,\varepsilon))+\lambda\mass{T}.
 \end{split}
\end{equation*}
This proves assertions (iii) and (iv) and hence the lemma.
\end{proof}
\subsection{Isoperimetric inequality for cycles of small mass}
In \cite{Wenger-GAFA} we have shown that every complete metric space $X$ which admits global cone type inequalities for $\intcurr_j(X)$,
$j=1,\dots,k$ admits an isoperimetric inequality of Euclidean type for all $T\in\intcurr_k(X)$ with $\bdry T=0$. (We point out that in this case
no restrictions are made on $T$ with respect to mass.) 
In this section we prove a similar result for spaces admitting $\delta$-cone type inequalities. These in general 
do not admit isoperimetric inequalities of Euclidean type for all cycles as is shown in \exref{example:counter-example-isoperimetric}.
The isoperimetric inequality of Euclidean type for cycles of small mass will enable us to apply \lemref{lemma:general-decomposition} in the proof of 
\thmref{theorem:flat-weak-convergence-higher-dimension}.
\bt\label{theorem:small-mass-isoperimetric-inequality}
 Let $(X,d)$ be a complete metric space admitting $\delta$-cone type inequalities for $\intcurr_j(X)$, $j=1,\dots,k$, for some $k\in\N$ and some $\delta>0$.
 Then there exist constants $C'>0$ and $D$ such that the following holds: If 
 $T\in\intcurr_k(X)$ satisfies $\bdry T=0$ and $\mass{T}\leq C'$ then there exists an $S\in\intcurr_{k+1}(X)$ such that $\bdry S=T$ and 
 \begin{equation*}
  \mass{S}\leq D[\mass{T}]^{\frac{k+1}{k}}.
 \end{equation*}
 The constant $D$ depends only on $k$ and on the constants of the $\delta$-cone inequalities, but not on $\delta$, whereas $C'$ also depends on $\delta$.
\et
The proof will show in particular that $C'=\infty$ in the case $\delta=\infty$ which returns Theorem 1.2 in \cite{Wenger-GAFA}.

\begin{example}\label{example:counter-example-isoperimetric}
 Let $X$ be the space obtained by gluing a $(k+1)$-dimensional
 half-sphere to $S^k\times[0,\infty)$ along $S^k\times\{0\}$ and the 
 equator of the half-sphere and endow it with the path metric coming from the
 path metrics on the sphere and on the cylinder. Then $X$ is ${\rm CAT}(1)$ but the cycle 
 $T:= \Lbrack S^k\times\{r\}\Rbrack\in\intcurr_{k}(X)$ has 
 \begin{equation*}
  \ifillnorm{T} = r\operatorname{Vol}(S^k) + \frac{1}{2}\operatorname{Vol}(S^{k+1}).
 \end{equation*}
\end{example}
\begin{proof}[{Proof of \thmref{theorem:small-mass-isoperimetric-inequality}}]
 We prove the theorem by induction on $k$. Let $k\geq 1$ and suppose in the case $k\geq 2$ that $X$ admits an isoperimetric inequality of Euclidean type
 with constant $C$ for cycles in $\intcurr_{k-1}(X)$ of mass no larger than $C''$. If $k=1$ then set $\lambda:=0$, $F:=1$, and $C':= \delta/4$. 
 If, on the other hand, $k\geq 2$ then define
 \begin{equation*}
  \lambda:= \min\left\{\frac{1}{6}, C\left(\frac{\omega_k}{2}\right)^{\frac{1}{k-1}}k^{\frac{k}{2k-2}}\right\}
 \end{equation*}
 as well as
 \begin{equation*}
  F:=\frac{\lambda^{k-1}}{C^{k-1}k^k}
 \end{equation*}
 and
 \begin{equation*}
  C':=\min\left\{F\delta^k4^{-k},C\lambda^{-1}(C'')^{\frac{k}{k-1}}\right\}.
 \end{equation*}
Let now $T\in\intcurr_k(X)$ satisfy $\bdry T=0$ and $\mass{T}<C'$.
Set  
\begin{equation*}
 r_0(y):=\max\{r\geq0: \|T\|(B(y,r))\geq Fr^k\}
\end{equation*}
and $r_1:= \frac{4}{3}r_0$. By \cite[Theorem 9]{Kirchheim} we have
\begin{equation*}
 \lim_{r\searrow 0}\frac{\hm^k(S_T\cap B(y,r))}{\omega_kr^k}=1
\end{equation*}
for $\hm^k$-almost all $y\in S_T$, the set $S_T$ being defined as in \eqref{equation:minimal-set}. 
Together with \thmref{theorem:mass-representation} and the 
fact that $F\leq\frac{\omega_k}{2k^{k/2}}$ this implies that the set of points $y\in X$ with $r_0(y)>0$ has full $\|T\|$-measure. 
It is then straight forward to check that $r_0$ and $r_1$ satisfy all the conditions of \lemref{lemma:general-decomposition}. 
We point out that with this particular choice of $r_0$ and $r_1$ conditions (a), (b), and (c) of \lemref{lemma:general-decomposition}  are
satisfied whenever $T\in\intcurr_k(X)$ is a cycle. If, moreover, $\mass{T}<C'$ then also (d) is fulfilled. With the choice of $C'$ made above the elements
$T_i$ in the decomposition of \lemref{lemma:general-decomposition} clearly satisfy
\begin{equation*}
 \diam(\spt T_i)\leq 4r_0(y_i)=4F^{-\frac{1}{k}}[\|T\|(B(y_i,r_0(y_i)))]^{\frac{1}{k}}\leq 4F^{-\frac{1}{k}}[\mass{T}]^{\frac{1}{k}}\leq\delta.
\end{equation*}
We successively apply \lemref{lemma:general-decomposition} to obtain (possibly finite)
sequences of cycles $(T_i)$, $(R_n)\subset\intcurr_k(X)$, of points $(y_i)\subset X$, and an increasing sequence $(N_n)\subset\N$
with the following properties:
\begin{itemize}
 \item $T= \sum_{i=1}^{N_n}T_i + R_n$
 \item $\diam(\spt T_i)\leq\delta$ and $\diam(\spt T_i)\leq E\mass{T_i}^{1/k}$ where $E:= 4[F(1-\lambda)]^{-\frac{1}{k}}$
 \item $\mass{R_n}\leq (1-\nu)^n\mass{T}$ where $\nu:= \left(1-\frac{1}{5^k}(1-\lambda)\right)$
 \item $\sum_{i=1}^\infty \mass{T_i}\leq \left[(1+\lambda)\sum_{i=0}^\infty(1-\nu)^i\right] \mass{T}=\frac{1+\lambda}{\nu}\mass{T}$.
\end{itemize}
An $S\in\intcurr_{k+1}(X)$ with $\bdry S=T$ and which satisfies the isoperimetric inequality is then constructed as follows. 
Since $\diam(\spt T_i)\leq \delta$ we can choose for each $T_i$ an $S_i\in\intcurr_{k+1}(X)$ with 
$\bdry S_i = T_i$ and such that
 \begin{equation}\label{eq:mass-S}
  \mass{S_i}\leq C_k\diam(\spt T_i)\mass{T_i} \leq C_kE\mass{T_i}^{\frac{k+1}{k}},
 \end{equation}
where $C_k$ denotes the constant of the $\delta$-cone inequality for $\intcurr_k(X)$.
The finiteness of $\sum_{i=1}^\infty\mass{T_i}$ implies that the sequence $S^n:= \sum_{i=1}^{N_n}S_i$
is a Cauchy sequence with respect to the mass norm because
 \begin{equation*}
  \mass{S^{n+q}-S^n}\leq C_kE\sum_{i=N_n+1}^\infty\mass{T_i}^{\frac{k+1}{k}}
                    \leq C_kE\left[\sum_{i=N_n+1}^\infty\mass{T_i}\right]^{\frac{k+1}{k}}
\;\;\overset{n\rightarrow\infty}{\longrightarrow}0.
 \end{equation*}
 Since $\intrectcurr_{k+1}(X)$ is closed in the Banach space $\curr_k(X)$ the sequence $S^n\in\intcurr_{k+1}(X)\subset\intrectcurr_{k+1}(X)$ 
 converges to a limit current $S\in\intrectcurr_{k+1}(X)$. As $T-\bdry S^n=R_n$ converges to $0$ it follows that
 $\bdry S=T$ and, in particular, that $S\in\intcurr_{k+1}(X)$.
 Finally, $S$ satisfies the isoperimetric inequality since
\begin{equation*}
 \mass{S}\leq \sum_{i=1}^\infty\mass{S_i}\leq C_kE\sum_{i=1}^\infty\mass{T_i}^{\frac{k+1}{k}}
         \leq C_kE\left(\frac{1+\lambda}{\nu}\right)^\frac{k+1}{k}\mass{T}^{\frac{k+1}{k}}.
\end{equation*}
This completes the proof.
\end{proof}
\subsection{Proof of \thmref{theorem:flat-weak-convergence-higher-dimension}}
The main step in the proof of \thmref{theorem:flat-weak-convergence-higher-dimension} will be to establish the following claim which clearly
implies the second statement of the theorem.\\*[1.5ex]
{{\bf Claim: }\it Let $1\leq k'\leq k$ and let $(T_m)\subset\intcurr_{k'}(X)$ be a bounded sequence satisfying $\bdry T_m=0$ for every $m\in\N$. 
If $(T_m)$ weakly converges to $0$ then $\ifillnorm{T_m}\to 0$ as $m\to\infty$.}\\*[1.5ex]
The above claim will be proved at the beginning of the proof of \thmref{theorem:flat-weak-convergence-higher-dimension} on page \pageref{claim-proof-page-number}.
We now need some preliminary results.
For this we define 
\begin{equation*}
 \varrho_T(\varepsilon):= \sup\{\|T\|(B(y,\varepsilon)):y\in X\}
\end{equation*}
whenever $T\in\intcurr_k(X)$ and $\varepsilon>0$.
\bp\label{proposition:fillvol-lambda-sequence}
 Let $(X,d)$ be a complete metric space and $\delta>0$. Let $k\in\N$ and suppose $X$ admits $\delta$-cone type inequalities 
 for $\intcurr_j(X)$, $j=1,\dots,k$.
 Then for every bounded sequence $(T_m)\subset\intcurr_k(X)$ of cycles and for every $\varepsilon>0$ the property
 \begin{equation*}
  \varrho_{T_m}(\varepsilon) \to 0\quad\text{as $m\to \infty$}
 \end{equation*}
 implies that
 \begin{equation*}
  \ifillnorm{T_m}\longrightarrow 0\quad\text{as $m\to\infty$.}
 \end{equation*}
\ep
\begin{proof}
 Let $\lambda_m>0$ be a sequence converging to $0$ and such that $\varrho_{T_m}(\varepsilon)\lambda_m^{1-k}\to 0$. 
 By \thmref{theorem:small-mass-isoperimetric-inequality} $X$ admits an isoperimetric inequality of Euclidean type with constant say $C$ for 
 cycles in $\intcurr_k(X)$ of mass no larger than $C'$.
 Define $F_m:=1$ if
 $k=1$ and 
 \begin{equation*}
  F_m:=\frac{\lambda_m^{k-1}}{C^{k-1}k^k}
 \end{equation*}
 in case $k\geq 2$. Denoting
 \begin{equation*}
  r_{0,{T_m}}(y):= \max\left\{r\in[0,\varepsilon]: \|T_m\|(B(y,r))\geq F_mr^k\right\}
 \end{equation*}
 one shows exactly as in the proof of \thmref{theorem:small-mass-isoperimetric-inequality} that for $m$ large enough $r_{0,{T_m}}>0$ almost everywhere. 
 Furthermore, by assumption,
 $r_{0,T_m}(y)\to 0$ uniformly in $y$ as $m\to\infty$. Set $r_{1,{T_m}}:=\frac{4}{3}r_{0,{T_m}}$ and note that for $m$ large enough 
 $T_m$, $r_{0,{T_m}}$ and $r_{1,{T_m}}$ satisfy all the
 conditions of \lemref{lemma:general-decomposition}. Therefore, there exist decompositions
 \begin{equation*}
  T_m=T_m^1+\dots+T_m^{N_m}+R_m
 \end{equation*}
 into the sum of cycles with the following properties:
 \begin{enumerate}
  \item $\max_{i=1,\dots,N_m}\diam(\spt T_m^i)\to 0$ as $m\to\infty$
  \item $\sum_{i=1}^{N_m}\mass{T_m^i}\leq (1+\lambda_m)\mass{T_m}$
  \item $\mass{R_m}\leq (1-\frac{1}{5^k}(1-\lambda_m))\mass{T_m}$
  \item $\|R_m\|(B(y,\varepsilon))\leq \|T\|(B(y,\varepsilon))+\lambda_m\mass{T_m}$.
 \end{enumerate}
 From (iv) we conclude, in particular, that $\varrho_{R_m}(\varepsilon)\leq \varrho_{T_m}(\varepsilon)+\lambda_m\mass{T_m}\to 0$ as $m\to\infty$.
 The proposition now follows by successive application of \lemref{lemma:general-decomposition} to $R_m$, together with 
 \thmref{theorem:small-mass-isoperimetric-inequality}.
\end{proof}
The next result is the higher dimensional analogue of the claim in the proof
of \thmref{theorem:0-dimension-flat-weak}.
\bp\label{proposition:reduction-density}
 Let $(X,d)$ and $\delta>0$ be as in \thmref{theorem:flat-weak-convergence-higher-dimension}. 
 Let $k\geq 1$ and suppose, in case $k\geq 2$, that the statement of the claim at the beginning of this section holds for $k':= k-1$. Let furthermore
 $(T_m)\subset\intcurr_k(X)$ be a bounded sequence of cycles weakly converging to $0$. 
 Then for every $\varepsilon>0$ small enough there exists a subsequence $(T_{m_j})$ and decompositions of $T_{m_j}$ into the
 sum of integral cycles
 \begin{equation*}
  T_{m_j}= U_j^1+\dots+U_j^{M_j}+V_j
 \end{equation*}
 satisfying the following properties:
 \begin{enumerate}
  \item $\ifillnorm{U_j^i}\leq \hat{C}\varepsilon\mass{U_j^i}$ for every $i\in\{1,\dots,M_j\}$ and $j\in\N$ and for a constant $\hat{C}$ depending only on the constants of 
        the $\delta$-cone type inequalities
  \item $\left[\sum_{i=1}^{M_j}(\mass{U_j^i} + \mass{V_j})\right] - \mass{T_{m_j}} \to 0$ as $j\to\infty$
  \item $\mass{V_j}\leq \mass{T_{m_j}} - \frac{1}{2}\varrho_{T_{m_j}}(\varepsilon/2)$ for every $j\in\N$
  \item $V_j$ weakly converges to $0$ as $j\to\infty$.
 \end{enumerate}
\ep
\begin{proof}
 \sloppy
 Let $(T_m)$ and $\varepsilon>0$ be given. In view of \thmref{theorem:small-mass-isoperimetric-inequality} and \propref{proposition:fillvol-lambda-sequence} 
 we may assume that $\liminf_{m\to\infty}\mass{T_m}>0$ and 
 $\limsup_{m\to\infty}\varrho_{T_m}(\varepsilon/2)>0$. We choose for every
 $m\in\N$ a point $y_m\in X$ with
 \begin{equation}\label{equation:lower-bound-mass}
  \|T_m\|(B(y_m,\varepsilon/2))\geq \frac{3}{4}\varrho_{T_m}(\varepsilon/2).
 \end{equation}
 In the following we distinguish two cases: First assume that after passing to a subsequence we have $y_m\in B(y_1,2\varepsilon)$ for every $m\geq 1$.
 By \lemref{lemma:boundary-slices-sequence} there exists $r\in(3\varepsilon,4\varepsilon)$ such that after passing yet to another subsequence we obtain that
 $T_m\rstr B(y_1,r)\in\intcurr_k(X)$ for all $m\in\N$, that $T_m\rstr B(y_1,r)\to 0$ as $m\to\infty$, and  
 $$\sup_{m\in\N}\mass{\bdry(T_m\rstr B(y_1,r))}<\infty.$$ By assumption, or by \thmref{theorem:0-dimension-flat-weak} in case $k=1$, there exist
 $S_m\in\intcurr_k(X)$ such that $\bdry S_m=\bdry(T_m\rstr B(y_1,r))$ and $\mass{S_m}\to 0$. By the remark following 
 \lemref{lemma:local-isoperimetric-support} (or the proof of \thmref{theorem:0-dimension-flat-weak} in case $k=1$) we may furthermore assume
 $\spt S_m\subset B(y_1,r+\varepsilon)$ for $m$ large enough. Clearly, $U_m^1:= T_m\rstr B(y_1,r)-S_m$ and $V_m:= T_m\rstr(X\ohne B(y_1,r))+S_m$ are integral cycles and
 $T_m=U_m^1+V_m$. Property (i) is then a consequence of the fact that $\spt U_m\subset B(y_1,r+\varepsilon)$ and of the $\delta$-cone type inequality whereas
 (ii) follows immediately from the fact that $\mass{S_m}\to 0$.
 Statement (iii) for $m$ large enough follows from \eqref{equation:lower-bound-mass} and from the fact that $\mass{S_m}\to 0$.
 Finally, (iv) holds because $T_m\rstr (X\ohne B(y_1,r))$ and $S_m$ both weakly converge to $0$. 
 This proves the proposition in the first case.\\
 We turn to the second case for
 which we assume no subsequence of $y_m$ with the property above exists. Then we may assume the balls $B(y_m,\varepsilon)$ to be pairwise disjoint and we set
 $A:= \{y_m: m\in\N\}$. By \lemref{lemma:boundary-slices-sequence} there exists an $r\in(\varepsilon/2,\varepsilon)$ such that after further passage to
 a subsequence $T_m\rstr A(r)\in\intcurr_k(X)$, $T_m\rstr A(r)\to 0$, and $\sup_m\mass{\bdry(T_m\rstr A(r))}<\infty.$
 As above there exist $S_m\in\intcurr_k(X)$ with $\bdry S_m=\bdry(T_m\rstr A(r))$ and $\mass{S_m}\to 0$. 
 Analogously, we may choose the $S_m$ in such a way that 
 $\spt S_m\subset A(\varepsilon')$ for some $r<\varepsilon'<\varepsilon$ for every
 $m$ large enough. It is then easy to see that $\bdry(T_m\rstr B(y_i,r))=\bdry(S_m\rstr B(y_i,\varepsilon))$ for every $i$. We set 
 $\tilde{U}^i_m:= T_m\rstr B(y_i,r) - S_m\rstr B(y_i,\varepsilon)$ for every $i\in\N$ and choose $M_m$ large enough such that 
 $U_m^{M_m}:=\sum_{i\geq M_m}\tilde{U}^i_m$ satisfies $\mass{U_m^{M_m}}\leq \varepsilon$. 
 Set $U_m^i:= \tilde{U}^i_m$ for $i<M_m$ and note that $U_m^i\in\intcurr_k(X)$ satisfies $\bdry U_m^i=0$ and, if $i<M_m$, 
 $\spt U_m^i\subset B(y_i,\varepsilon)$. 
 Then we can use the $\delta$-cone type inequality in case $i<M_m$, 
 respectively \thmref{theorem:small-mass-isoperimetric-inequality} in case $i=M_m$, to conclude that the $U_m^i$ indeed satisfy property (i). 
 Set $V_m:= T_m\rstr (X\ohne A(r))+S_m$.
 It is now easy to check that also properties (ii), (iii), and (iv) hold. Indeed, we have
 \begin{equation*}
  \begin{split}
   \sum_{i=1}^{M_m}&\mass{U_m^i} + \mass{V_m}\\
                   &\leq \sum_{i=1}^{\infty}\|T_m\|(B(y_i,r)) + \mass{S_m} + \|T_m\|(X\ohne A(r)) + \mass{S_m}\\
                   &= \mass{T_m}+2\mass{S_m}
  \end{split}
 \end{equation*}
 from which (ii) follows. As for (iii) we use \eqref{equation:lower-bound-mass} to obtain
 \begin{equation*}
  \begin{split}
   \mass{V_m}&\leq \|T_m\|(X\ohne A(r))+\mass{S_m}\\
             &=\mass{T_m}-\|T_m\|(A(r)) +\mass{S_m}\\
             &\leq \mass{T_m}-\frac{3}{4}\varrho_{T_m}(\varepsilon/2) +\mass{S_m}.
  \end{split}
 \end{equation*}
Since $\mass{S_m}\to 0$ and $\limsup_{m\to\infty}\varrho_{T_m}(\varepsilon/2)>0$ statement (iii) follows.
Finally, property (iv) is a consequence of the fact that $T_m\rstr A(r)$ and therefore also  $T_m\rstr (X\ohne A(r))$ weakly converges to $0$ 
 and so does $S_m$. 
\end{proof}
We now put everything together to prove \thmref{theorem:flat-weak-convergence-higher-dimension}. The main task will be to prove the claim made at the beginning
of this section. From this the theorem will follow quite easily.
\begin{proof}[Proof of \thmref{theorem:flat-weak-convergence-higher-dimension}]
It is clear that convergence with respect to the flat distance implies weak convergence. In order to prove the converse we prove the 
claim\label{claim-proof-page-number} stated 
at the beginning of this section. This will clearly establish the second statement of the theorem.
In order to prove the claim we first note that it is enough to find a subsequence $T_{m_j}$ satisfying $\ifillnorm{T_{m_j}}\to0$. 
Let $\varepsilon>0$ be small enough. The proof is on induction on $k'$. The arguments for the case $k'=1$ and for the induction step are 
almost identical. Let therefore $1\leq k'\leq k$ and suppose, in the case $k'\geq 2$, that the claim holds for $k'-1$.
Successive application of \propref{proposition:reduction-density} together with an argument involving diagonal sequences yields a subsequence $(T_{m_j})$ and 
decompositions $T_{m_j}= U_j^1+\dots+U_j^{N_j}+R_j$ into the sum of integral cycles satisfying the following properties:
 \begin{enumerate}
  \item $\sum_{i=1}^{N_j}\mass{U_j^i} + \mass{R_j} \leq 2\mass{T_{m_j}}$ for every $j\in\N$
  \item $\ifillnorm{U_j^i}\leq \hat{C}\varepsilon\mass{U_j^i}$ for every $i\in\{1,\dots,N_j\}$ and $j\in\N$ for a constant $\hat{C}<\infty$
  \item $\varrho_{R_j}(\varepsilon/2)\to 0$ as $j\to\infty$.
 \end{enumerate}
 Using the $\delta$-cone type inequality, respectively \propref{proposition:fillvol-lambda-sequence} for $R_j$, we easily see that 
 $$\liminf_{j\to\infty}\ifillnorm{T_{m_j}}\leq \tilde{C}\varepsilon$$ for some constant $\tilde{C}$ depending only on $\sup\mass{T_m}$ and the constants of the 
 $\delta$-cone type inequalities. Since $\varepsilon$ was arbitrary this proves the claim and hence the second statement of the theorem.\\
 It remains to prove the first statement. Let $(T_m)\subset\intcurr_k(X)$ be a bounded sequence converging weakly to some 
 $T\in\intcurr_k(X)$. Set $T'_m:= \bdry T_m-\bdry T$ and observe that $T'_m\in\intcurr_{k-1}(X)$, $\bdry T'_m=0$, and $\sup_{m\in\N}\mass{T'_m}<\infty$. 
 Since $T'_m$ weakly converges to $0$ the claim with $k'=k-1$ (or \thmref{theorem:0-dimension-flat-weak} in case $k=1$) implies the existence of 
 $R_m\in\intcurr_k(X)$ with $\bdry R_m = T'_m$ and $\mass{R_m}\to 0$. 
 Hence the $T''_m:= T_m-T-R_m\in\intcurr_k(X)$ form a bounded sequence of cycles and by the claim there exist $S_m\in\intcurr_{k+1}(X)$ with 
 $\bdry S_m=T''_m$ 
 and $\mass{S_m}\to 0$. Since $T_m-T=R_m+\bdry S_m$ it follows that
\begin{equation*}
 \iflnorm{T_m-T}\leq \mass{R_m}+\mass{S_m}\to 0.
\end{equation*}
 This concludes the proof of the theorem.
\end{proof}

\section{Applications}\label{section:applications}
\fussy
From Theorems \ref{theorem:0-dimension-flat-weak} and \ref{theorem:flat-weak-convergence-higher-dimension} we infer the useful result below on 
absolutely area minimizing currents.
\bt\label{theorem:minimal-weak}
 Let $k\geq 1$ and let $(X,d)$ be a metric space satisfying the conditions of \thmref{theorem:flat-weak-convergence-higher-dimension}. 
 If $(S_m)\subset\intcurr_k(X)$ is a bounded sequence converging weakly to $S\in\intcurr_k(X)$ with
 \begin{equation*}
  \mass{S_m} - \ifillnorm{\bdry S_m}\to 0
 \end{equation*}
 then $S$ is absolutely area minimizing. In particular we have $$\mass{S}=\liminf_{n\to\infty}\mass{S_n}.$$
\et
\begin{proof}
 By Theorem \ref{theorem:0-dimension-flat-weak} and \ref{theorem:flat-weak-convergence-higher-dimension} for every $m$ there exists  $R_m\in\intcurr_k(X)$ 
 with $\bdry R_m=\bdry S-\bdry S_m$ and such that $\mass{R_m}\to 0$ as $m\to\infty$. 
 For a current $V\in\intcurr_k(X)$ with $\bdry V= \bdry S$ we see that $\bdry(V-R_m)=\bdry S_m$ and hence
 \begin{equation*}
   \mass{V} \geq \mass{V-R_m} - \mass{R_m}
    \geq \ifillnorm{\bdry S_m} - \mass{R_m}
 \end{equation*}
 from which we conclude 
 \begin{equation*}
  \mass{V}\geq \liminf_{m\to\infty}\mass{S_m}\geq \mass{S}.
 \end{equation*}
\end{proof}
A consequence of Theorems \ref{theorem:0-dimension-flat-weak} and \ref{theorem:flat-weak-convergence-higher-dimension} and the compactness theorem for integral 
currents \cite[Theorems 5.2 and 8.5]{Ambr-Kirch-curr} is the following result.
\bt\label{theorem:compactness-flat-topology}
 Let $k\geq 0$ be an integer and $(X,d)$ a quasiconvex compact metric space which, in case $k\geq 1$, admits local cone type 
inequalities for $\intcurr_j(X)$, $j=1,\dots,k$. Then $\{T\in\intcurr_k(X): \nmass{T}\leq L\}$ is compact with respect to $d_{\mathcal F}$ for every $L<\infty$.
\et
Finally, using the equivalence of flat and weak convergence, we prove the existence of minimal elements in fixed Lipschitz homology classes in compact
quasiconvex metric spaces admitting local cone type inequalities, as stated in \thmref{theorem:minimal-homology}. 
As for the definitions we set for a complete metric space $(X,d)$ and $A\subset X$ closed
\begin{equation*}
 {\mathcal Z}_k(X,A):= \{T\in\intcurr_k(X):\spt(\bdry T)\subset A, \spt T\text{ compact}\}
\end{equation*}
\begin{equation*}
 {\mathcal B}_k(X,A):= \{R+\bdry S: R\in\intcurr_k(A), S\in\intcurr_{k+1}(X), \spt S, \spt R\text{ compact}\}.
\end{equation*}
It is then easy to check (see \cite[4.4.1 (1)--(5), (7)]{Federer} and \cite[p.~10]{Eilenberg-Steenrod} for the excision property) that
\begin{equation*}
 {\bf H}_k(X,A):= {\mathcal Z}_k(X,A)/{\mathcal B}_k(X,A)
\end{equation*}
satisfies the Eilenberg-Steenrod axioms for a homology. See \cite[4.4.5]{Federer} for the corresponding definitions in case $X\subset \R^n$ is a compact local
Lipschitz neighborhood retract. 
\begin{proof}[{Proof of \thmref{theorem:minimal-homology}}]
 Let $T_m\in{\mathcal Z}_k(X,A)$ be a minimizing sequence subject to the condition $[T_m]=c$ for all $m\in\N$. 
 By the compactness theorem for integral currents there exists 
 a subsequence converging weakly to some $T\in\intcurr_k(X)$. It is clear that $\spt(\bdry T)\subset A$ and hence $T\in{\mathcal Z}_k(X,A)$. 
 In order to check that $[T]=c$ we view $\bdry T_m$ and $\bdry T$ as elements of $\intcurr_{k-1}(A)$ and apply 
 \thmref{theorem:flat-weak-convergence-higher-dimension} to conclude the existence of $R_m\in\intcurr_k(A)$ such that
 $\bdry T - \bdry T_m=\bdry R_m$ for all $m\in\N$ and $\mass{R_m}\to 0$ as $m\to\infty$. Using \thmref{theorem:flat-weak-convergence-higher-dimension} 
 one more time we obtain a sequence $S_m\in\intcurr_{k+1}(X)$ satisfying $\bdry S_m=T-T_m-R_m$ for 
 all $m\in\N$ and $\mass{S_m}\to 0$ as $m\to\infty$. Thus we have $T-T_m=R_m+\bdry S_m$ with $\spt R_m\subset A$ which proves that indeed 
 $T-T_m\in{\mathcal B}_k(X,A)$ and hence $[T]=[T_m]=c$.
\end{proof}


\begin{thebibliography}{GG7}
   \bibitem[AK]{Ambr-Kirch-curr} L.~Ambrosio, B.~Kirchheim:
    {\it Currents in metric spaces}, Acta Math. 185 (2000), no. 1, 1--80.
  \bibitem[Ba]{Ballmann} W.~Ballmann: {\it Lectures on spaces of nonpositive curvature}, DMV Seminar Band 25, Birkh\"{a}user, 1995.
  \bibitem[BBI]{Burago-Burago-Ivanov}D.~Burago, Yu.~Burago, S.~Ivanov: {\it A course in metric geometry}, Graduate Studies in Mathematics, Vol. 33,
   Amer. Math. Soc., Providence, Rhode Island, 2001.
  \bibitem[BrH]{Bridson-Haefliger} M.~R.~Bridson, A.~Haefliger 
    {\it Metric Spaces of Non-Positive Curvature}, 
    Grundlehren der mathematischen Wissenschaften 319, Springer, 1999.
  \bibitem[DG]{deGiorgi} E.~De Giorgi: {\it Problema di Plateau generale e funzionali geodetici}, Atti Sem Mat. Fis. Univ. Modena, 43 (1995), 282--292.
  \bibitem[ES]{Eilenberg-Steenrod}S.~Eilenberg, N.~Steenrod: {\it Foundations of algebraic topology}, Princeton University Press, 1952.
  \bibitem[Ek]{Ekeland} I.~Ekeland: {\it On the variational principle}, 
    J.~Math.~Anal.~Appl. 47 (1974), 324--353.
  \bibitem[Fe]{Federer} H.~Federer: {\it Geometric Measure Theory}, Springer 1969, 1996.
  \bibitem[FF]{Fed-Flem} H.~Federer, W.~H.~Fleming:
    {\it Normal and integral currents}, Ann. Math. 72 (1960), 458--520. 
  \bibitem[Ki]{Kirchheim} B.~Kirchheim: {\it Rectifiable metric spaces: local structure and regularity
         of the Hausdorff measure}, Proc. Am. Math. Soc. 121 (1994), no. 1, 113--123.
  \bibitem[We]{Wenger-GAFA} S.~Wenger: {\it Isoperimetric inequalities of Euclidean type in metric spaces}, accepted for publication in GAFA.
 \end{thebibliography}
\end{document}